\documentclass[12pt]{amsart}

\usepackage{esint}
\usepackage{stackrel}

\usepackage{pb-diagram}

\usepackage{amsfonts}
\usepackage{amsmath}
\usepackage{amssymb}

\newcommand{\dbar}{\overline{\partial}}

\newcommand{\ddt}[1]{\frac{\partial #1}{\partial t}}
\newcommand{\R}{\mathcal{R}}
\newcommand{\Sc}{\mathcal{S}}

\usepackage{xcolor}

\newcommand{\A}{\mathcal{A}}

\newcommand{\ddbar}{\sqrt{-1}\partial\dbar}

\newtheorem{theorem}{Theorem}[section]
\newtheorem{proposition}{Proposition}[section]
\newtheorem{lemma}{Lemma}[section]

\newtheorem{example}{Example}[section]
\newtheorem{definition}{Definition}[section]

\newtheorem{remark}{Remark}[section]

\def \C {\textbf C}

\def \R {\mathcal R}

\def \U {\mathcal U}

\def \O {\mathcal O}

\begin{document}

\title{Uniqueness of tangent cone of K\"ahler-Einstein metrics on singular varieties with crepant singularities }

\author{Xin Fu}

\address{Department of Mathematics, University of California, Irvine, CA 92617}

\email{fux6@uci.edu}


\begin{abstract} Let $(X,L)$ be a polarized Calabi-Yau variety (or canonical polarized variety)   with crepant singularity. Suppose $\omega_{KE}\in c_1(L)$ (or $\omega_{KE}\in c_1(K_X))$ is the unique Ricci flat current (or K\"aher-Einstein current with negative scalar curvature) with local bounded potential constructed in \cite{EGZ}, we show that the local tangent at any point $p\in X$ of metric $\omega_{KE}$  is unique. 
\end{abstract}

\maketitle

\section{Introduction}
Uniqueness of tangent cone problems have been of great interest in many geometric contexts. For example, minimal varieties, Yang-Mills connections and Einstein metrics. The influential paper of Simon \cite{Simon} provides a framework to deal with the uniqueness problem given one of the tangent cone has smooth link and leads to a proof of uniqueness of tangent cone of minimal variety if one tangent cone has smooth link. In this note, we study the tangent cone of K\"ahler Einstein metric.
Let $(X_i,g_i,p_i)$ be a sequence noncollapsed Riemannian manifold with Ricci bounded below, by Gromov compactness Theorem, we have
$$(X_i,g_i,p_i)\stackrel{GH}{\longrightarrow}(X_\infty,g_\infty,p_\infty),$$
where $X_\infty$ is a metric length space. 
The regularity of $X_\infty$ is intensively studied in the last thirty years mainly by Cheeger, Colding and Naber \cite{CC1,CC2,CC3,Colding,CJN,LiuSze1}. If the sequence of Riemannian manifolds $X_i$ have two sided Ricci bound, the limit space $X_\infty$ has more regularity(cf. \cite{CCT,DS,CN,JN}). One important tool of studying the regularity  of $X_\infty$ is to study its tangent cone.  Give a singular point $p\in X_\infty$ and a sequence of positive constants $r_i\to 0$, it is shown that the rescaled sequence $(X_\infty, r^{-2}_ig_\infty, p_\infty)$ converge by sequence to a metric cone \cite{CC1}. Based on the tangent cone, there is regular-singular decomposition of $X_\infty=\mathcal R\cup\mathcal S$, where $\mathcal R$ consists of points admitting a tangent cone which is isometric to $\mathbb R^{2n}$ and $\mathcal S$ is the complement of $\mathcal R$. In the K\"ahler-Einstein context, $\mathcal R$ also has a complex structure induced from $X_i$ and $\mathcal S$ has Hausdorff codimension at least $4$ by Cheeger-Colding-Tian.  Now a natural question is whether the tangent cone is independent of the rescaling sequence $r_i$. It is known from Perelman \cite{Pere} that, with  Ricci only bounded below, the tangent cone is not unique. On the other hand, if $X_i$ are Einstein manifolds with bounded Einstein constant, the uniqueness of tangent cone at infinity problem is firstly studied in the  work of Cheeger-Tian \cite{CT} and was confirmed under the assumption that one of the tangent cone has isolated singularity in Colding-Minicozzi \cite{CM} by using a Lojasiewicz-Simon type argument. Generic uniqueness is also obtained in the recent work Cheeger-Jiang-Naber \cite{CJN}. We should remark that there are  few results confirming the uniqueness of tangent cone without assuming that one of the tangent cone has smooth link. One exception is: if we further restrict  to the case that $(X_i,g_i)$ are polarized K\"ahler-Einstein manifolds, Donaldson-Sun \cite{DS2} proves the uniqueness of tangent cone in full generality by using an completely new algebraic geometric approach. 

In this short note, we want to study the tangent cone of K\"ahler-Einstein metrics on varieties with crepant singularties.  More precisely, suppose $(X,L)$ is a polarized Calabi-Yau (or canonical polarized) variety (indeed allowing a more general class of singularity), by  Eyssidieux-Guedj-Zeriah \cite{EGZ} (see also \cite{Z}), there is a K\"ahler Einstein current with bounded local potential $\omega_{KE}\in c_1(L)$. 
Such K\"ahler Einstein current is unique, thus canonical,  and its associated K\"ahler-Einstein metric is smooth on the regular
part of the underlying variety. The same existence result is also obtained
in \cite{ZYG,To}. Ricci flat metrics on varieties with crepant singularities are also intensively studied as metric degeneration of smooth Calabi Yau manifolds and conifold transitions, we refer the interested reader to the work \cite{To1,RUZ,RZ,S1,S,CPY}.  Going back to the K\"ahler current $\omega_{KE}$, it is natural to study the geometry of $\omega_{KE}$. If we denote the metric completion of $(X_{reg},\omega_{KE})$ by $X_\infty$, it is
shown by Rong-Zhang \cite{RZ} that the metric regular part $\mathcal R$ of $X_\infty$ coincides with the complex analytic regular part $X_{reg}$.
It is further shown in Song \cite{S} that, the metric completion $X_\infty$ of $(X_{reg},\omega_{KE})$ is indeed homeomorphic to the underlying complex variety $X$. For such an canonical K\"ahler-Einstein metric on $X$, we shall study  its tangent cone structure. Fix a singular point $p\in X_\infty$, the main goal of this paper is to show that the local tangent cone at $p$ is unique. We  remark that the main difference of our result with \cite{DS2} is that there, the K\"ahler-Einstein metric on the singular variety is approximated by a sequence of smooth K\"ahler-Einstein  manifold coming from polarization and here, in this note, it is not necessary true that such an approximation can be achieved for $(X,\omega_{KE})$. Now let us state our main result.
\begin{theorem}\label{main1} Let $(X,L)$ be a polarized Calabi-Yau variety with crepant singularity and let $\omega_{KE}\in c_1(L)$ be the Ricci flat metric constructed in \cite{EGZ}, then for any point $p\in X$, the local tangent cone of metric $\omega_{KE}$ at $p$ is unique.
\end{theorem}
We have a similar result on canonical polarized variety with crepant singularity.
\begin{theorem}\label{main22}
 Let $X$ be a canonical polarized  variety with crepant singularity and let $\omega_{KE}\in K_X$ be the negative K\"ahler-Einstein metric constructed in \cite{EGZ} with locally bounded potential, then for any point $p\in X$, the tangent cone of metric $\omega_{KE}$ at $p$ is unique.
\end{theorem}
\begin{remark}
In Theorem \ref{main1} and Theorem \ref{main22},  The variety $X$ is shown to be homeomorphic to the metric completion $X_\infty$ of $(X_{reg},\omega_{KE})$ in by Song in \cite{S} and moreover, the metric regular part $\mathcal R$ of $X_{reg}$ is exactly the complex analytic regular part of the underlying variety $X$. Therefore, there is no difference between a complex analytic singular point of $X$ and a metric singular point of $X_\infty$.
\end{remark}
\section{Preliminaries}

In this section, we will recall some basic definitions and notations . 

\begin{definition} Let $X$ be a normal projective variety $X$. If $K_X$ is a Cartier $\mathbb{Q}$-divisor and if $\pi: \tilde X \rightarrow X$ is a resolution of $X$, then there exist $a_i \in \mathbb{Q}$ with 
$$K_{\tilde X} = \pi^* K_X + \sum_i a_i E_i, $$
where $E_i$ ranges over all exceptional prime divisors of $\pi$. 
$X$ is said to have terminal singularities (canonical singularities,  log terminal singularities, log canonical), if all $a_i >0$ ( $a_i \geq 0$, $a_i > -1 $, $a_i\geq -1$).  In particular, $X$ is said to have crepant singularities if $$K_{\tilde X} = \pi^* K_X, $$ and $\pi$ is called a crepant resolution of $X$.

\end{definition}

\begin{example} All surface $A$-$D$-$E$ singularities are crepant singularities. 
\end{example}




%
%
%
%
%
%
%
%
%

Now we  define projective Calabi-Yau varieties as our main interest.

\begin{definition}\label{a} A projective Calabi-Yau variety $X$ is a normal projective variety with canonical singularities and torsion  canonical divisor $K_X$, i.e.  $qK_X$ is a trivial line bundle for some positive integer $q$
.\end{definition}

By Kawamata's base point free theorem  \cite{KMM}, the following  theorem is well known. 

\begin{lemma} 

Let $L\rightarrow X$ be a holomorphic line bundle over a projective Calabi-Yau variety $X$. If $L$ is big and nef, then it is   semi-ample.

\end{lemma}
Immediately, we have many Cababi-Yau varieties with crepant singularities through the construction below.

\begin{example} 

Suppose $L$ is a big and nef line bundle over a projective Calabi-Yau manifold $X$. Then the linear system $|L^k|$ induces a unique surjective birational morphism 
$$\Phi: X \rightarrow Y$$
for sufficiently large $k\in \mathcal{F}(X,L)$ such that $Y$ is a projective Calabi-Yau variety with crepant singularities.

\end{example}

All projective Calabi-Yau varieties with conifold singularities  admit a crepant resolution, although crepant resolutions are not necessarily unique and they are related by flops. 



\section{Analytic tools}
In this section, we collect various crucial estimates on the singular Calabi-Yau varieties with crepant singularities from Song \cite{S}, which inspires our work. And we do not claim any originality for these results.
\subsection{Construction of K\"ahler-Einstein metrics by perturbation}

In this subsection, we review how K\"ahler-Einstein metrics are constructed. More precisely, in this subsection, we will assume that $X$ is an $n$-dimensional  normal projective Calabi-Yau variety with crepant singularities and $L \rightarrow X$ is an ample $\mathbb{Q}$-line bundle over $X$.   We also let $$\pi: X' \rightarrow X $$ be a crepant resolution of $X$ and $L' = \pi^* L$. We also denote the smooth and singular part of $X$ by $X_{reg}$ and $X_{sing}$. In particular, $X_{sing}$ is an analytic subvariety of $X$ of complex codimension at least $2$. 

It is shown in \cite{EGZ} that there exists a unique Ricci-flat K\"ahler current  in $c_1(L)$ with bounded local potential. We review  the construction. Let  $\A$ be an arbitrary ample line bundle on $X'$. Since $X'$ is a Calabi-Yau manifold, there exists  a smooth volume form $\Omega'$ on $X'$satisfying 
$$\ddbar \log \Omega' = 0, ~ \int_{X'} \Omega' = [L']^n.$$ Obviously, $\Omega' = \pi^* \Omega$ for some smooth Calabi-Yau volume form on $X'$.
We can choose $\chi' =\pi^*\chi \in c_1(L')$ be the Fubini-Study metric induced by an projective embedding from the base point linear system of  $L^k$ for some sufficiently large k. 
We also choose $\omega_{\A} \in [\A] $ a fixed smooth K\"ahler metric on $X'$. We then consider the following Monge-Ampere equation
\begin{equation}\label{CY}
(\chi + e^{-t} \omega_{\A} + \ddbar \varphi_t)^n = e^{c_t} \Omega, ~\int_{X'}  \varphi_t dg_t=0~, t\in [0, \infty),  
\end{equation}
where $g_t$ is the K\"ahler metric associated to the K\"ahler form $\omega_t=\chi+ e^{-t} \omega_{\A} + \ddbar \varphi_t$,  $c_t$ is normalization constant satisfying $ e^{c_t} [L']^n =[L'+e^{-t} \A]^n$. It is straightforward to see that  $c_t= O(e^{-t}  )$ for sufficiently large  $t$. 

By Yau's  solution \cite{Y1} to the Calabi conjecture, equation (\ref{CY}) admits a smooth solution $\varphi_t$ for $t$ large enough and $g_t$ is a smooth Ricci-flat K\"ahler metric. Moreover, if we denote $E:=\pi^{-1}(X_{sing})$ to be the exceptional locus of a crepant resolution, then for any $k>0$ and $K\subset\subset X'\setminus E$, there exists constant $C_{K,k}$ such that for all $t\in [0,\infty)$, we have $\|\varphi_t\|_{C^k(K)}\leq C_{K,k}$. 
%

%
%
%

Therefore, by letting $t\rightarrow \infty$,  $\varphi_t$ will converge locally and uniformly to a function $\varphi_\infty \in PSH(X, \chi)\cap L^\infty(X)  \cap C^\infty(X_{reg})$ solving the equation 
\begin{equation}\label{SMA}
(\chi + \ddbar \varphi_\infty)^n =\Omega.
\end{equation}

%
%
%
%
The K\"ahler current $\omega_\infty=\chi + \ddbar \varphi_\infty$ lies in class $c_1(L)$ with bounded local potential. Here and after, we denote $g_\infty$ to be the  K\"ahler metric induced by $\omega_\infty$, which is smooth on $X_{reg}$.


\subsection{$L^2$-estimates}

In this subsection, we collect various $L^2$-estimates for global sections in $H^0(X, L^k)$.

Without confusion, we can identify the  quantities on $(X_{reg}, L)$ and $(X'\setminus E, L')$ through the resolution map $\pi$.
Let $h_{FS}$ be the smooth fixed hermitian metric on $L$ induced from some embedding of $L^k$ for some $k$ large such that $Ric(h_{FS}) = \chi$. We define the singular hermitian metric $h_\infty$ and $h_\infty'$ on $L$ and $L'$ by
$$h_\infty = e^{-\varphi_\infty} h_{FS}, ~ h_\infty' = \pi^*h_\infty.$$
Easily, we have $Ric (h_\infty) = \omega_\infty$. 


%
%
%


We define the rescaled norm $|| \cdot ||_{L^{\infty, \sharp}}$ and $|| \cdot ||_{L^{2, \sharp}}$ for $s\in H^0(X, L^k)$ with respect to the hermitian metric $h_\infty^k$ and $kg_\infty$. In the following lemmas, we derive several $L^2$ type estimates with respect to the rescaled norm. There estimates will be crucial in Proposition \ref{C0}. The first lemma says for holomorphic sections of line bundle $L$, its $L^\infty$ norm can be controlled by its $L^2$ norm.

\begin{lemma} \label{l2}\cite[Proposition 3.3]{S}  There exists  $K >0$ such that if  $s\in H^0(X, L^k)$ for $k\geq 1$, then 
\begin{equation}\label{l1}
\|s\|_{L^{\infty, \sharp}} \leq K \|s\|_{L^{2, \sharp}}.
\end{equation}

\end{lemma}

The following lemma gives a pointwise bound for the gradient of $s\in H^0(X, L^k)$.

\begin{lemma}\cite[Proposition 3.4]{S} \label{grad} There exists  $K >0$ such that if  $s\in H^0(X, L^k)$ for $k\geq 1$, then 
  $$ \|\nabla s\|_{L^{\infty, \sharp}} \leq K \|s\|_{L^{2, \sharp}} $$

\end{lemma}

We will also need the following version of $L^2$-estimates due to Demailly  for a big and net line bundle (singular Hermitian metric) over a weakly pseudoconvex K\"ahler manifold. 

\begin{lemma}\cite[Corollary 5.3]{De} \label{demailly} Let $X$ be a $n$-dimensional smooth complex manifold (possibly open) equipped with a  complete smooth K\"ahler metric $\theta$ and $\omega$ be another  K\"ahler metric on $X$. Also let $L$ be a holomorphic line bundle over $X$ equipped with a possibly singular hermitian metric $h$ such that $Ric(h) = -\ddbar \log h \geq \delta \omega$ in current sense for some $\delta>0$. Then for every  $L$-valued $(n,1)$-form $\tau$ satisfying 
$$ \dbar \tau =0, ~ \int_X |\tau|^2_{h, \omega} ~\omega^n<\infty,$$
where $|\tau|_{h, \omega}^2 = tr_{\omega} \left( \frac{ h \tau \overline\tau}{\omega^n} \right)$, 
there exists an  $L$-valued $(n, 0)$-form $u$ such that $\dbar u = \tau$ and 
\begin{equation}
\int_X |u|^2_h ~\omega^n \leq \frac{1}{2\pi \delta} \int_X |\tau|_{h, \omega}^2 ~\omega^n.
\end{equation}

\end{lemma}

We need a generalization of standard Hormander's $L^2$ estimate on manifold to the singular varieties of our interest. In a recent interesting work of Liu-Szekelyhidi \cite{LiSze2}, $L^2$ type estimate is obtained on singular tangent cone and leads to deep structure result for the tangent cone. Here we aim to solve $\dbar$-equation on $X$ itself.

\begin{proposition}\cite[Propsition 3.5]{S} \label{L2Sing} Let $X$ be a projective Calabi-Yau variety with crepant singularities. If $L$ is an ample line bundle over $X$ equipped with a hermitian metric such that $\omega=Ric(h) \in c_1(L) $ is the unique Ricci-flat K\"ahler current on $X$ with bounded local potentials, then for any smooth $L$-valued $(0,1)$-form $\tau$ satisfying

\begin{enumerate}

\item $\dbar \tau =0$,  

\item $Supp ~\tau \subset \subset X_{reg}$, 

\end{enumerate}
there exists an $L$-valued section $u$ such that $\dbar u = \tau$ and $$ \int_X |u|^2_h ~\omega^n \leq \frac{1}{2\pi} \int_X |\tau|^2_h~ \omega^n. $$

\end{proposition}

\begin{proof} Since this proposition is crucial to us, we copy the proof here for convenience. Let $\pi: X' \rightarrow X$ be the crepant resolution of $X$.   Since $L$ is ample on $X$, by Kodaira's lemma, there exists a divisor $D$ on $X'$ such that $L' - \epsilon [D]$ is ample for all $\epsilon>0$, where $L'= \pi^*L$.    Let $s_D$ be the defining section of $D$ and $h_D$ a smooth hermitian metric on the line bundle induced by $[D]$ satisfying $$\chi + \epsilon \ddbar \log h_D >0$$ 
for all sufficiently small $\epsilon>0$, where $\chi \in c_1(L')$ is a smooth closed semi-positive $(1,1)$-form as the pullback of the Fubini-Study metric from the linear system of $|(L')^k|$ for some sufficiently large $k$. We consider the following Monge-Ampere equation 
$$(\chi + \epsilon \ddbar \log h_D + \ddbar \varphi_\epsilon)^n = c_\epsilon \Omega, ~ \int_{X'} \varphi_\epsilon \Omega =0 $$
where $\Omega$ is a fixed smooth Calabi-Yau volume form on $X'$ and $c_\epsilon \int_{X'}\Omega = [L - \epsilon D]^n \rightarrow [L^n]$ as $\epsilon \rightarrow 0$. 
Obviously  $$\omega_\epsilon = \chi + \epsilon \ddbar \log h_D + \ddbar \varphi_\epsilon $$ is the unique Ricci flat K\"ahler metric in $[L'-\epsilon D]$. $\varphi_\epsilon$, $h_\epsilon$ and $\omega_\epsilon$ converges to $\varphi$, $h=h_{FS}e^{-\varphi}$ and $\omega=\chi + \ddbar \varphi$ weakly globally on $X$ and smoothly on $X'\setminus D$ as $\epsilon \rightarrow 0$.

We can identify $s\in H^0(X', L')$ and $s \in H^0(X', (L' -K_{X'})+  K_{X'})$ and the hermitian metrics on $L'$ are also hermitian metrics on $L'-K_{X'}$ because $K_{X'}$ is numerically trivial.

We now define $h_\epsilon = h_{FS} e^{ - \epsilon \log |s_D|_{h_D}^2 - \varphi_\epsilon}$ to be a singular hermitian metric on $L'$, in particular, 
$$ Ric(h_\epsilon)= -\ddbar \log h_{\epsilon} = \omega_\epsilon + [s_D] \geq \omega_\epsilon$$
in the sense of currents.   Hence we can apply Lemma \ref{demailly} to $X', L', \omega_\epsilon, h_\epsilon, \tau' = \pi^*\tau $ (Here $\tau'$ is identified as a $L'-K_{X'}$ valued $(n,1)$ form). Note that  $\tau$ is smooth and 
$$\lim_{\epsilon\rightarrow 0} \int_{X'} |\tau'|_{h_\epsilon}^2 \omega_\epsilon^n =\int_{X'} |\tau'|_{h}^2 \omega^n < \infty$$
since $\tau'$ vanishes in a neighborhood of the exceptional locus of $\pi$ and so $h_\epsilon$ and $\omega_\epsilon$ converges smoothly on the support of $\tau'$.
By Lemma \ref{demailly}, there exists $u_\epsilon $ on $X'$ such that 
$$\dbar u_\epsilon = \tau', ~~ \int_{X'} |u_{\epsilon}|^2_{h_\epsilon} \omega_\epsilon^n \leq \frac{1}{2\pi} \int_{X'} |\tau'|_{h_\epsilon}^2 \omega_\epsilon^n. $$

Also $h$ and $h_{FS}$ are uniformly equivalent since $\varphi$ is uniformly bounded, $$\int_{X'} |u_\epsilon|^2_h~ \omega^n \leq C \int_{X'} |u_\epsilon|^2_{h_{FS}} \Omega \leq C \int_{X'} |u_\epsilon|^2_{h_\epsilon} ~\omega_\epsilon^n $$
is uniformly bounded for all $\epsilon>0$. 

Hence we can take a subsequence of $u_{\epsilon_v} $ converging weakly in $L^2(X', \Omega)$ to $u \in L^2(X, \Omega)$ as $\epsilon_v \rightarrow 0$ and 
$$\dbar u =\dbar \lim_{\epsilon \rightarrow 0} u_{\epsilon_v} = \tau'$$
on $X'$ in the distribution sense.
On the other hand, $u$ is bounded uniformly in $L^2(X, h_{FS}e^{-\epsilon \log |s_D|^2_{h_D}} \Omega)$ for any $\epsilon>0$ and $\omega_\epsilon$ converges to $\omega$ in $C^\infty(X\setminus D)$,  hence  
$$\int_{X'} |u|^2_h \omega^n \leq \int_{X'} |\tau' |_h^2  \omega^n.  $$
The proof is complete after pushing $u$ to $X$.
\end{proof}


\medskip


By the gradient estimate Lemma \ref{grad} and the fact that $X_\infty$ is a complete metric space, any holomorphic section $\tau\in H^0(X, L^k)$ can be continuously extended from $X_{reg}$ to $X_\infty$. Therefore we can use $H^0(X, L^k)$ to construct maps $\Phi_{k,\sigma}$ from $X_\infty $ to projective spaces as follows: 
$$\Phi_{k, \sigma} : z\in X_\infty \rightarrow  [ \sigma_1(z), ..., \sigma_{d_k+1}(z) ] \in\mathbb{CP}^{d_k} $$
where $\{\sigma_1, ..., \sigma_{d_k+1} \}$ is a basis of $H^0(X, L^k)$.

\begin{remark}
Since $L$ is semi-ample, $\Phi_{k, \sigma}$ is stabilized for sufficiently large and divisible. 
\end{remark}
\subsection{Gromov-Hausdorff limits}

In this subsection, we can summarize the results of Song \cite{S}  in the following theorem.  By previous discussions, we can construct a map $\Phi_{k,\sigma}: X_\infty\to\mathbb{CP}^{d_k}$, the following theorem essentially says that for $k$ large, $\Phi_{k,\sigma}$ is a homeomorphism  from $X_\infty$ onto $X$.

\begin{theorem}\cite[Theorem 1.1]{S} \label{GHlim} Let $g_t$ be the unique Ricci-flat K\"ahler metric solving equation (\ref{CY}) on $X'$. Then $(X', g_t)$ converges in Gromov-Hausdorff topology to a unique metric length space $(X_\infty, d_\infty)$. Furthermore,  

\smallskip

\begin{enumerate}

\item $(X_\infty, d_\infty)$ is the metric completion of $(X_{reg}, g_\infty)$.

\smallskip

\item $X_\infty = \mathcal{R}  \cup \mathcal{S} $, where $ \R$ is the regular part  of $(X_\infty, d_\infty)$ and $\Sc$ is the singular set. In particular, $\Sc$ is closed and its Hausdorff dimension is no greater than $n-4$. 

\smallskip

\item $\R$  is  convex in $(X_\infty, d_\infty)$ and 
%
%
$\R= X_{reg}.$ 
%
\item $\Phi=\Phi_{k,\sigma}:X_\infty\to X$ is a homeomorpism for $k$ sufficiently large.
\end{enumerate}
\end{theorem}


\medskip





\section{uniqueness of tangent cone: Calabi-Yau case}
The  goal of this section is to prove Theorem \ref{main1}. We will follow the arguments of Donaldson-Sun \cite{DS2} very closely. On the other hand, we should remark that, in \cite{DS2}, they proved that if $X_\infty$ is the Gromov-Hausdorff limit of a sequence of noncollapsed polarized K\"ahler-Einstein manifold, then the tangent cone is unique. In particular, the tangent cone can also    be approximated by a sequence of noncollapsed K\"ahler-Einstein manifold and the smooth approximation gives many useful analytic tools to study the tangent cone. In our case, our key observation is that, instead of using the smooth polarized K\"ahler Einstein manifolds to approximate  the tangent cone, we directly use the rescaling of $X_\infty$  to approximate the tangent cone. Therefore, various uniform analytic estimates must be derived directly on the singular metric space $(X_\infty,g_\infty)$ and its rescaled space. We will point out where new analytic estimates, when compared with Donaldson-Sun, are derived in the process of proof.
\begin{theorem}(Restate Theorem \ref{main1}) Let $X_\infty$ be the Gromov-Hausdorff limit in Theorem \ref{GHlim}. Fix any point $p\in X_{\infty}$, then tangent cone $C(Y_p)$ at $p$ is an normal affine variety and  moreover tangent cone is unique both as a metric cone and as a complex variety.
\end{theorem}

By Cheeger-Colding theory, we know that for any sequence of positive integers $r_i\to 0$, \begin{equation}\label{data}(X_i:=X_\infty, p_i:=p, g_i:=r_i^{-2}g_\infty) \stackrel{GH}{\longrightarrow}(C(Y),g_{C(Y)},o).\end{equation}
Here $g_i$ is the unique Ricci flat metric in class $c_1(L_i)$, where $L_i:=L^{r^{-1}_i}$. Our goal is to show that the limit space is independent of the choice of $\{r_i\}_{i=1}^{i=\infty}$. 
 It will be convenient to introduce the following space consisting of pointed metric space.
 
 Fix a constant $\lambda=1/\sqrt{2}$, 
 we define the tangent cone space at $p$ as:  
 \begin{equation}\label{cone}\mathcal C_p:=\{ C(Y)|C(Y) \,\textnormal{is the pointed Gromov-Hausdorff limits of} \,(X_i,g_i, p_i)\, \textnormal{for some}\, r_i\to 0\}\end{equation}
 Let us also introduce space $\mathcal T_p$.
 \begin{equation}\label{cone1}
 \mathcal T_p:=\mathcal C_p\cup\{(X_i,g_i,p_i)|(X_i,g_i,p_i)=(X_\infty,\lambda^{-2i},p)|\,i\in\mathbb Z\}
 \end{equation}
The following remark says that to exhaust all tangent cones, it is enough to consider the rescaled sequences with rescaling factor $\lambda^{-2i}, i\in\mathbb Z$.
 \begin{remark} \label{T}
Let $C(Y)$ and $C(Y')$ be two tangent cones at $p$ defined by two sequences of positive numbers $\{c_k\}, \{d_k\}$ respectively. If for some constant $A>0$ and $k$ sufficiently large, we have $A^{-1} \leq c_k/d_k\leq A$,  then $C(Y)$ and $C(Y')$ are isometric.
\end{remark}
 \begin{lemma}\label{codim4}
  Recall that we have defined $(X_i,g_i,p_i):=(X_\infty, r^{-2}_ig_\infty, p)$ for some $r_i>0$. If  $$(X_i,g_i, p_i)\stackrel{GH}{\longrightarrow}(C(Y),g_{C(Y)},o)$$
in the Gromov-Hausdorff sense, then $g_i$ converge smoothly to the metric regular part $\mathcal R$ of $C(Y)$ and the metric singular part $\mathcal S$ of $C(Y)$ has Hausdorff codimension at least 4.
 \end{lemma}
 \begin{proof}
  The $C^\infty$ convergence of metric to the regular part is a direct consequence of Anderson's regularity Theorem. By the K\"ahler condition $\nabla_{g_i} J_i=0$, the complex structure also converge smoothly. Also since $X_\infty$ is Gromov-Hausdorff limit of a sequence of  noncollapsed  Einstein manifold by Theorem \ref{GHlim} (But not polarized Einstein manifold), therefore its tangent cone is also the Gromov-Hausdorff limit of a sequence of noncollapsed  Einstein manifolds, then by Cheeger-Naber \cite{CN}, $C(Y)$ has codimension 4 singularity.  
 \end{proof}
 The following lemma will be needed in several places. 
\begin{lemma}\label{snotr}
Suppose we are in the context of Lemma \ref{codim4}. If $q_i$ is a sequence of points in the metric singular part of $X_i$ converging to $q$ in $C(Y)$, then $q$ is  also in the metric singular part of $C(Y)$.\end{lemma}
\begin{proof}
If not, then $q$ lies in the smooth part of $C(Y)$. Then for any $\epsilon>0$, there exists a sufficiently small $r>0$ such that $\frac{Vol(B_r(q))}{Vol(B^{\mathbb E}_r(0))}>1-\frac{\epsilon}{2}$. Then by Colding's volume continuity theorem, we have   $\frac{Vol(B_r(q))}{Vol(B^{\mathbb E}_r(0))}>1-{\epsilon}$ for $i$ large enough. Recall that $X_\infty$ itself is the Gromov-Hausdorff limit of a sequence of noncollapsed Einstein manifolds. Then by Anderson's gap theorem \cite{Anderson}, $q_i$ must be a smooth point in the metric sense. Contradiction. 
\end{proof}
\begin{remark}
Lemma \ref{codim4} and Lemma \ref{snotr} also hold for iterated tangent cones i.e. tangent cones of $C(Y)$.
\end{remark}
With the above lemma, we can graft objects supported on the smooth part of the tangent cone to the smooth part of the approximating spaces. More precisely,  for any $j$, if we  
denote the complement of the $j^{-1}$ neighborhood of the
singular set in $B(o,j)\subset C(Y)$ by $U_j$, then there is a  embedding $\chi_i$ of $U_j$ into the smooth part of $B(p_i,j)$ for $i$ sufficiently large. Moreover, we have
for all $x\in U_j, d(x, \chi_i(x)) < j^{-1}$  and
$\|\chi^*_ig_i-g_{C(Y)}\|
_{C^j} +\|\chi_i^*J_i-J_{C(Y)}\|_{
C^j} < j^{-1},$ where the norm is taken with respect to $g_{C(Y)}$.

\subsection{Partial $C^0$ estimate} Then following partial $C^0$ type \cite{T,DS} is very useful. 
We introduce the $H$-property introduced by Donaldson-Sun \cite{DS}. 

\begin{definition} We consider the follow data $(p_*, D, U, J, L , g, h, A)$ satisfying 

\begin{enumerate}
\item  $(p_*, U, J, g)$  is an open bounded K\"ahler manifold with a complex structure $J$, a K\"ahler metric $g$ and a base point $p_* \in D\subset \subset U$ for an open set $D$, 

\item $L \rightarrow U$ is a hermitian line bundle equipped with a hermitian metric $h$ and $A$ is the connection induced by the hermitian metric $h$ on $L$, with  its curvature $\Omega(A) = g$.

\end{enumerate}

\noindent The data $(p_*, D, U, L, J, g, h, A)$ is said to satisfy the  $H$-condition if there exist $C>0$ and  a compactly supported smooth section $\sigma: U\rightarrow L$   satisfying 

\begin{enumerate}

\item [$H_1$:]   $\|\sigma \|_{L^2(U)} < (2\pi)^{n/2}$,

\item [$H_2$:]  $|\sigma(p_*) | >3/4$,

\item [$H_3$:] for any holomorphic section $\tau$ of $L$ over a neighborhood of $\overline{D}$,
$$|\tau(p_*)| \leq C (\| \dbar \tau \|_{L^{2n+1}(D)} + ||\tau||_{L^2(D)} ),$$

\item [$H_4$:]  $\|\dbar \sigma \|_{L^2(U)} < \min \left( \frac{1}{8\sqrt{2} C},  10^{-20} \right)$,

\item [$H_5$:]  $||\dbar \sigma ||_{L^{2n+1}(D)} \leq \frac{1}{8C}$.

\end{enumerate}

\end{definition}

Here all the norms are taken with respect to $h$ and $g$. The constant $C$ in the $H$-condition depends on the choice $(p_*, D, U, J, L, g, h)$.

An important fact about the $H$-condition is that it can be achieved on the tangent cone when the singular set of tangent cone has Harsdorff codimension 4. 

\begin{lemma} \cite{DS}\label{Hcone} Let $q \in C(Y_{reg})$. If $3/4 <e^{ -dist(q,o)^2}  <1$, then for any $\epsilon>0$, there exists $U\subset \subset C(Y_{reg})\setminus \{o\}$ and an open neighborhood $D\subset\subset U$ of $q$ such that 
 $(q, D, U, L_C, J_C, g_C, h_C, A_C)$ satisfies the $H$-condition.

\end{lemma}

The following proposition from \cite{DS} establishes the stability of the $H$-condition for perturbation of the curvature and the complex structure. 

\begin{proposition}  \label{Hstability} Suppose $(p_*, D, U, J_C, L_C, g_C, h_C, A_C )$ constructed as above in Lemma \ref{Hcone} satisfies the $H$-condition. There exist $\epsilon>0$ and $m \in \mathbb{Z^+}$ such that for any collection of data $(p_*, D, U, J, g, h, A)$ if 
$$|| g - g_C ||_{C^0(U(m)) } + ||J-J_C ||_{C^0(U(m)) } < \epsilon, $$
then for some $1\leq l \leq m$, 
$$(p_*, D, U, \mu_l^*J, \mu_l^*L, \mu_l^*g, \mu_l^*h, \mu_l^*A)  $$
satisfies the $H$-condition.

\end{proposition}
Let us recall the converging Gromov-Hausdorff sequence considered  in (\ref{data}), i.e. we have a sequence of positive integers $r_i\to 0$ and a sequence of metric spaces
$$(X_i:=X_\infty, p_i:=p, g_i:=r_i^{-2}g_\infty, L_i:=L^{r_i})$$
converging to a tangent cone $(C(Y),g_{C(Y)},o)$.
We need the following key local partial  $C^0$ estimate essentially proved by Donald-Sun  in \cite{DS2}.  

 \begin{proposition}\label{C0} For any point $p\in X_\infty$, there are constants $\lambda_1, \lambda_2\in (0, 1)$ with $\lambda_1>\lambda_2$ such that if some rescaled sequence of  metric space $(X_i, g_i, p_i)$ as above converge to a tangent cone $(C(Y),o)$ in the Gromov-Hausdorff sense, then there are constants $k$, $C$, $N$, $l_1, \cdots, l_N$  with the following effects: 
 
\begin{enumerate}
\item[(1)] For $i$ sufficiently large, there is a $s_i\in H^0(X_i,L_i^k)$  such that,  under the rescaled norm $k\omega_i$, $|s_i(x)|\geq 3/4$ when $d_i(p, x)\leq \lambda_1$, and $||s_i||_{L^{2}}\leq (2\pi)^{\frac{n}{2}}$.
\item[(2)]  There are holomorphic sections $\sigma_{j, i}\in H^0(X_i,L_i^{kl_j})$ for $j=1,\cdots, N$. 
\item[(3)] The map $$F_i:X_i\to\mathbb C^N, F(x)=\Big(\sigma_{1, i}/s_i^{l_1}(x),\cdots,\sigma_{N, i}/s_i^{l_N}(x)\Big) $$ satisfies $|F_i(x)|> 1/3$ when $d_i(p, x)=\lambda_1$, and $|F_i(x)|\leq 1/50$ when $d_i(p, x)\leq \lambda_2$. Here $|\cdot|$ standards for the Euclidean norm.  
\item[(4)] $|\nabla F_i|\leq C$ on smooth points of $B_{\lambda_1}$, and Lipschitz with Lipschitz constant $C$. 

\end{enumerate}
\end{proposition}
\begin{proof} If  $(X_i,p_i,g_i,L_i)$ is a sequence of smooth noncollapsed polarized K\"ahler-Einstein manifolds, the result is proved in \cite{DS2}.  We point out that the proof of \cite{DS2} relies on two important key ingredients: the first one is the existence of  Sasakian Einstein tangent cone  for the limit metric space of $(X_i,p_i,g_i,L_i)$ and its codimension 4 regularity structure.  The second one is several versions of uniform $L^2$ type estimates for polarized K\"ahler-Einstein metrics.  
Now we sketch the proof,  consider one of the tangent cone $C(Z)$ at point $o$ (Here $C(Z)$ is an iterated tangent cone). By Lemma \ref{Hcone}, the $H$-condition data $(p_*,D,U,J,L,g,h,A)$ can be constructed on the cone $C(Z)$. Let us also fix a sequence of embedding $\chi_i$ from $U$ to the regular part of suitable rescaling of $X_i$ (cf. Lemma \ref{codim4}) realizing the Gromov-Hausdorff distance.  Then by the stability of $H$-condition  under small smooth perturbation (Lemma \ref{Hstability}), we may assume, for $i$ large,  $(\chi_i^{-1})^{*}(p_*,D,U,J,L,g,h,A)$ also satisfy the $H$-condition.  Now by uniform  $L^2$ estimate (Lemma \ref{l2}, Lemma \ref{grad} and Proposition \ref{L2Sing}), following \cite{DS},   there is some $k\leq k(o)$ such that for all $i>i(o)$, if $x$ is a point in $X_i$ with
$d(x, o)\leq \lambda_1$, then there exists a section $s_i$ of line bundle $L_i^{k}$ satisfying $|s_i| (x)>b(o).$ For such sections $s_i$, it is usually named peak sections. To finish the construction of $F_i$,  it suffices to construct finite peak sections $\sigma_{j}$ with center landing  near the boundary $\partial B_{\lambda_1}(o)$ and graft them back to $\sigma_{j,i}$ on $X_i$ (cf. \cite[Proposition 16]{CDS2}). Then it is easy to verify items $(1)$ to $(4)$ are satisfied.
\end{proof}
\begin{remark}
We remark here that item (4) in the above lemma use the gradient estimate for holomorphic sections (cf. Lemma \ref{grad}) and the metric completeness of $X_\infty$.
\end{remark}


There are some direct consequences of Proposition \ref{C0}.

Let $B_{1/4}$ be a Euclidean ball with radius $\frac{1}{4}$ in $\mathbb C^N$. Let
 $\Omega_i:=F_i^{-1}(B_{1/4})$ and $W_i:=F_i(\Omega_i)$,  then by item (3), the fiber $F_i^{-1}(x)$ is a compact analytic set for any $x\in B$. By item (1) the ample line bundle $L_i^{ka_i}$ is trivial over $F_i^{-1}(B)$, so $F_i^{-1}(x)$ can have at most finite points. Combined with Remmert proper mapping theorem \cite[Theorem 8.8]{Demailly}, $F_i$ is a finite map from $\Omega_i$ onto an analytic set $W_i$.  By the gradient estimate item (4), volume of $W_i$ induced  by the Euclidean norm  on $\mathbb C^N$ are uniformly bounded.  By Bishop's compactness theorem for complex analytic sets, $W_i$ converges by sequence to a limit analytic set $W$. Item (4) also implies that we can take a limit of $F_i$ and obtain a Lipschitz map $F$ from an open neighborhood $\Omega$ of $p$ onto $W$. 
 
 
 
 

 Now let us define a sheaf structure on the open neighbourhood $\Omega$ of $p$. We first define a presheaf $\mathcal F$ on $\Omega$ as follows: for any open set $\mathcal U\subset \Omega$ and $\mathcal U_i\subset\Omega_i$ converging to $\mathcal U$ , define
 \begin{equation}\label{sheaf}\mathcal F(\mathcal U):=\{f|f=\lim_i f_i,\,\,\,\textnormal{where} \,\,\,f_i\in\mathcal O(\mathcal U_i) \}
 \end{equation}
 Then let $\mathcal O$ be the sheafification of $\mathcal F$.
 We first derive uniform $L^\infty$ and gradient estimate for functions $f\in\mathcal O$ defined as above.
 
 \begin{proposition} \label{gradholo} Recall that the space $\mathcal T_p$ consists of rescaled metric spaces of $X_\infty$ and all tangent cones of $X_\infty$ at $p$ (cf. equation \ref{cone}). 
Let $(Y,p)$ be a pointed metric space in $\mathcal T_p$ and $B_r(p)$ be a metric ball with radius $r$ contained in $Y$.  Then there are  constants $C_0=C_0(n,r)$, $C_1=C_1(n,r)$, such that for any holomorphic function $f\in\mathcal O(B_r(p))$, we have  
$$|f(p)|\leq C_0 |f|_{L^2(B_r(p))}, $$
$$|\nabla f(p)|\leq C_1|f|_{L^2(B_r(p))}.$$
\end{proposition}
\begin{remark}
At this moment, when $Y$ is a tangent cone, then $f\in\mathcal O(B_r(p))$ is understood as in the sense of \ref{sheaf}.
\end{remark}

\begin{proof} For simplicity, we may assume $r=1$. Also notice that we only need to prove the theorem on all rescaled space $(X_\infty, r_i^{-2}g_\infty, p)$, because if we further pass to limit, the tangent cone case will also be proved.
 Now suppose $f$ is a holomorphic function on the unit ball near $p$. Since $(X_\infty,g_\infty,p)$ is the Gromov-Hausdorff limit of noncollapsed  manifold $(X',g_t),$  with diameter bounded above and Ricci curvature bounded below (cf. \cite[Proposition 3.1]{S}), so the Sobolev constant is uniformly controlled in terms of $t$ and $i$. (Note that rescaling will make the Sobolev constant better). 
 Now fix a resolution map $\pi:X'\to X$ and also a sequence of domain $(B_t,g_t)$ converging to $B$ as in Theorem \ref{GHlim}. Now we pull function $f$ back to $X'$. Then standard Moser iteration gives
 $$|f(p)|\leq C_1|f|_{L^2(B_t,g_t)}.$$
Now since $f$ is bounded and we can find open neighbourhoods $V$ of the singular set $\mathcal S$ with $Vol(V)$ as small as we want, then by Colding's volume continuity theorem, the $L^2$ norm of $|f|$ on the singular set can be made arbitrary small. So we have $\lim_{t\to\infty}|f|_{L^2(B_t,g_t)}=|f|_{L^2(B,g_\infty)}$. This proves the first inequality. If $q\in\mathcal R$ is a point on the one quarter ball, then gradient estimate of Cheng-Yau \cite{CY} for positive harmonic function implies that 
$$|\nabla_{g_t} f(q)|\leq C_1|f|_{L^\infty(B_{\frac{1}{2}})}\leq C_1|f|_{L^2(B)}.$$
 Let $t\to\infty$, we have $|\nabla f(q)|\leq C_1|f|_{L^2(B)}$. Now by the metric completeness of $X_\infty$ and the density of the regular set $\mathcal R$ in $X_\infty$, inequality  $|\nabla f(p)|\leq C_1|f|_{L^2(B)}$ holds in the Lipschitz sense for singular point $p\in X_\infty$.
 
\end{proof}

By adding finite functions  $\hat F\in \mathcal O(\Omega)$ as new components to obtain a new map $F'=(F, \hat F):  \Omega\rightarrow W'\subset \mathbb C^{N+1},$ Donaldson-Sun proves that  $\Omega$ is homeomorphic to its image through map $F'$ . 
 \begin{lemma}\label{homeo}\cite[Proposition 2.3]{DS2}
By further shrinking $\Omega$ if necessary, $F$ is a homeomorphism from $\Omega$ to an analytic variety $W$. Moreover, $F$ maps the metric regular set $\mathcal R$ to the smooth part of $W$.
\end{lemma}

Next we want to show that $W$ is normal. To do this,
 we need to derive a local Hormander $L^2$ type estimate on the domain $\Omega_i$. 
 \begin{proposition}\label{L2loc}
 Let $(\Omega_i, L^{r^{-1}_i}, \omega_i:=r_i^{-2}\omega_\infty)$ be defined as above, then for any smooth $L^i$-valued $(0,1)$-form $\tau$ satisfying 
\begin{enumerate}

\item $\dbar \tau =0$,  

\item $Supp ~\tau \subset \subset X_{reg}\cap\Omega_i$, 

\end{enumerate}
there exists an $L$-valued section $u$ such that $\dbar u = \tau$ and $$ \int_{\Omega_i} |u|^2_h ~\omega_i^n \leq \frac{1}{2\pi} \int_{\Omega_i} |\tau|^2_h~ \omega_i^n. $$

 \end{proposition}
 \begin{proof}
 The proof is the same as Proposition \ref{L2Sing}. Firstly, notice that, we have a resolution map $\pi: X'\to X$ and also a holomorphic map $F_i:\Omega_i\to\mathbb C^N$, hence $F_i\circ\pi $ is a holomorphic map form $\pi^{-1}\Omega_i$ to $\mathbb C^N$. So $\varphi:=(F_i\circ \pi)^*|z|^2$ will be a smooth PSH function on $\pi^{-1}\Omega_i$. Since $\Omega_i:=\{|z|^2<\frac{1}{16}\}$, $\pi^{-1}\Omega_i:=\{\varphi<\frac{1}{16}\}$. Note that $\ddbar\varphi$ is only a smooth semipositive (1-1) form on $\pi^{-1}(\Omega_i)$, but we can still solve $\bar\partial$-equation with $L^2$ estimate  on the weakly pseudoconvex K\"ahler manifold $\pi^{-1}(\Omega)$ by Lemma \ref{demailly}. Therefore we can use the same argument as Proposition  \ref{L2Sing} to complete the proof i.e. solving $\dbar$ on $\pi^{-1}\Omega_i$ with a family of singular metric $\omega_\epsilon$ (already constructed globally in Proposition \ref{L2Sing}) and then limit the family of sections $u_\epsilon$ to $\Omega_i$. 
 \end{proof}
 With the local $L^2$ type estimate in hand. Given a bounded holomorphic function $f$ on the regular part $\mathcal R$ of $W$, We can first pull it back to domain $\Omega_i$ with a good cutoff function, then solve $\bar\partial$-equation with controlled estimate to get a holomorphic function $f_i$ on $\Omega_i$ and at last limit $f_i$ to some function $\tilde f$ on $\Omega$, which extends $f$ across the singular set. So we have,
 \begin{lemma}\label{normal}\cite[Propsition 2.4]{DS2} $W$ is a normal variety and $F$ induce a sheaf isomorphism: $F^*:\mathcal O_W\to\mathcal O_\Omega$.
 \end{lemma}

\subsection{More Analytic tools}
Before we proceed to study the tangent cone, we collect several useful analytic consequences of partial $C^0$ Proposition \ref{C0}. We remark that these results will be fundamental which relates  the tangent cone with its Gromov-Hausdorff approximation.

\begin{lemma}
 There are constants $\lambda_1, \lambda_2\in (0,1)$ with $\lambda_1>\lambda_2$ such that if $(Y,p)$ is a pointed metric space in $\mathcal T_p$, then there is
an open set $U$ in $Y$  containing the closure of the ball $B_{\lambda_1}(p)$, and a
holomorphic map $F: U\to\mathbb C^N$ satisfying:
\begin{enumerate}
\item $|F(x)|>\frac{1}{3}$ when $d(p,x)=\lambda_1$,
\item $|F(x)|<\frac{1}{50}$ when $d(p,x)<\lambda_2$,
\item $|\nabla F(x)|<C$ when $d(p,x)<\lambda_1$.
\end{enumerate}
\end{lemma}
\begin{proof}
We will focus on the case $(Y,p)$ is not a tangent cone and the tangent cone case will follow by taking limit. Now suppose there are no such uniform constants $\lambda_1,\lambda_2$. Then for a sequence of positive number
$\lambda_{2,i}<\lambda_{1,i}\to 0$, there are a sequence of  pointed space $(Y_i,p_i)$ which do not admit a map $F$ satisfying the property stated in the lemma. By compactness and Proposition \ref{C0}, this is contradiction. \end{proof}
\begin{lemma} \label{Approximatinglemma} There are constants $\lambda_1, \lambda_2\in (0,1)$ with $\lambda_1>\lambda_2$ such that if $(C(Y),o)$ is a tangent cone and
  $f$ is a holomorphic function defined on $B_{\lambda_1}(o)$. Then for any sequence of pointed metric spaces $(Y_i,p_i)$ in $\mathcal T_p$ which converge to $C(Y)$ in the Gromov-Hausdorff topology,   there is a holomorphic function $f_i$ defined on  $B_{\lambda_2}(p_i)$, such that  $f_i$ converges to  $f$ uniformly over $B_{\lambda_2}(p)$.
  
\end{lemma}
\begin{proof} This lemma is essentially due to  \cite[Lemma 2.9]{DS2}.  We give more details here. Without lose of generality, we also assume $\lambda_1=1$. Recall the definition of $\mathcal T_p$ in (\ref{cone1}). We will focus on the case $Y_i$ are all tangent cones.  Suppose we have a sequence of tangent cone $(C(Y_i),o_i)$ converge to $(C(Y),o)$ in the Gromov-Hausdorff sense. Then for any sequence of unit ball $(B_i,o_i)\subset C(Y_i)\to (B,o)\subset C(Y)$,
 let us fix a  sequence $B_{i,j}\subset (X_i,r_j^{-2}g_\infty)\to B_i$ as $j=1,\cdots,n$ realizing the Gromov-Hausdorff distance. \\

\textbf{Claim:} For any $\eta>0$ and $f\in\mathcal O(B)$, there exists $\epsilon(\eta)$ such that if $d_{GH}(\tilde B,B)<\epsilon$ with $\tilde B\subset X_\infty$ (In particular, this means that the ball $\tilde B$ is a neighbourhood of $p$ in $X_\infty$ rather than some subset of a tangent cone), then exist $\tilde f\in\mathcal O(\tilde B)$ ($\tilde f$ to be constructed will only be defined on an open set contained in $\tilde B$ and containing $\tilde B_{\lambda_2}$, which will be enough for your purpose. For simplicity, we assume $\tilde f$ is defined on $\tilde B$.) such that $|\tilde f-f|\leq\eta$ on $B_{\lambda_2}$. Here we assume a Gromov-Hausdorff approximation map is choosen and $\tilde f$ is pulled back to $B_{\lambda_2}.$\\

Assuming this claim, for any $\eta>0$, when $i$ is sufficiently large, there exists an integer $n_i$ such that when $j\geq n_i$, there are functions $f_{i,j}\in \mathcal O(B_{i,j})$ with $|f_{i,j}-f|\leq\eta$ on $B_{\lambda_2}$. By the gradient estimate Lemma \ref{gradholo}, we can limit $f_{i,j}$ to find $f_i\in\mathcal O(B_i)$ with $|f_i-f|\leq\eta$ on $B_{\lambda_2}$. Taking $\eta=\frac{1}{i}$, we find functions $f_i\in\mathcal O(B_i)$ satisfying $|f_i-f|\leq\frac{1}{i}$ on $B_{\lambda_2}$. The  lemma will follow by limiting $f_i$.\\

\textbf{Proof of claim:} We first fix a sequence $\epsilon_i\to 0.$ Using the fact that
the singular set $\mathcal S$ has Hausdorff codimension at least 4, as in \cite{DS}, we
have a sequence of cut-off functions $\beta_i$ on $\Omega$ satisfying $\beta_i$ is
supported in the complement of a neighborhood of $\mathcal S$ and  $\beta_i=1$ outside
the $\epsilon_i$ neighborhood of $\mathcal S$. Moreover $\|\nabla \beta_i\|_{L^2}\leq \epsilon_i$. Let us recall the context of Proposition \ref{C0}, where we have constructed a sequence of domain $\Omega_i$ containing $B_{\lambda_1}$ and a sequence of Gromov-Hausdorff approximation map $\chi_i:\Omega\to\Omega_i$.  For the given holomorphic function $f$, we use the maps $\chi_i$
to pull $f\beta_i$ back to $\Omega_i$ and construct a smooth section $\sigma_i = (\chi_i^{-1})^*
(f\beta_i)s_i$
of $L^{kr_i^{-1}}
$ over $\Omega_i$ with $\|\bar\partial\sigma_i\|_{L^2}\to 0$.  By the local $L^2$ estimate Proposition \ref{L2loc}, we can find a section $\tau_i$ solving $\bar\partial\tau_i=\bar\partial\sigma_i$ with the following estimate $\|\tau_i\|_{L^2}\leq\|\bar\partial\sigma_i\|_{L^2}$. Moreover,   on any compact set with fixed distance away from the singular set, $\bar\partial\tau_i=0$ when $i$ is sufficiently large. Hence $|\nabla\tau_i|\to 0$
on these compact sets. 

Now we define the desired function $f_i:=\frac{\sigma_i-\tau_i}{s_i}.$
Then $f_i$ is a sequence of holomorphic functions with uniformly bounded $L^2$ norm. Then by the gradient estimate Lemma \ref{gradholo}, $f_i\to f$ pointwisely on $\Omega$ and uniformly on $B_{\lambda_2}$, which finish the proof.


\end{proof}


\begin{lemma} \cite[Proposition 2.12]{DS2} \label{Fgeneric}
Let $(C(Y),o)$ be a tangent cone and $F: B_{\lambda_1}(o)\rightarrow \mathbb C^N$ be a holomorphic embedding. Then for any pointed metric space $(Y_i,p_i)\in\mathcal T_p$ converging to $(C(Y),o)$, there is a sequence of generic one-to one holomorphic map $F_i: B_{\lambda_1}(p_i)\to\mathbb C^N$  such that the image  $F_i(B_{\lambda_2}(p_i))$ converges to $F(B_{\lambda_2}(p))$ as local complex analytic sets in $\mathbb C^N$. 
\end{lemma}

Now we study the convergence property of holomorphic functions. Let $B_i\subset X_i$ be a sequence of ball converging to $B\subset C(Y)$ and $f_i\in\mathcal O(B_i)$ be a sequence of holomorphic functions converging to $f\in\mathcal O(B)$. 
 By the estimate in Lemma \ref{gradholo},   the convergence is uniform over any compact subset of $B$ if the $L^2$ norm of $f_i$ is uniformly bounded.  We always have the following weak convergence of $L^2$ norm, 
 $$||f||_{L^2(B)}\leq \liminf_{i\rightarrow\infty} ||f_i||_{L^2(B_i)}.$$
Also if
 $$||f||_{L^2(B)}= \lim_{i\rightarrow\infty} ||f_i||_{L^2(B_i)},$$
 we say $f_i$ converges strongly to $f$.

\begin{lemma}\cite[Lemma 2.16]{DS2}\label{L2Converge} 
\begin{enumerate}
\item[(1)] If $f_i$ converges uniformly to $f$, then $f_i$ converges strongly to $f$. \item[(2)]
  If $g_i$ converges weakly to $g$ and $f_i$ converges strongly to $f$.  If $f$ also extends to a holomorphic    function over a larger ball , then $\int_{B} f\bar g=\lim_{i\rightarrow\infty}\int_{B_i} f_i\bar g_i. $
\end{enumerate}
\end{lemma}
\begin{proof}
The proof will be identical to \cite[Lemma 2.16]{DS2} except that there, for a smooth domain $V\subset B\cap\mathcal R$, they use a smooth sequence of domain $V_i$ to approximate $V$. Here a priori, $V_i$ is not necessary smooth. By Lemma \ref{snotr}, $V_i$ is indeed smooth.
\end{proof}
 To finish this subsection, we show one more property of $F$.
 \begin{proposition}\label{STOS} F maps the metric singular set $\mathcal S$ to complex analytic singular set of $W$.
 \end{proposition}
  Again, if $\Omega_i$ is smooth, this lemma is proved in \cite{DS2}. It seems that the technique there could not be applied in our set-up. Instead, we use an idea from Liu-Szekelyhidi \cite{LiuSze1}. We record the following lemma to characterize metric smoothness, which will be used to  prove Proposition \ref{STOS}.

\begin{lemma} \label{criteria}\cite[Proposition 5.2]{LiuSze1} Assume that $p\in C(Y)$ is not a metric regular point. Then there
exist $\epsilon > 0$ and $r_1 > 0$ such that for all $r < r_1$, if nonzero holomorphic
functions $f_1, \cdots, f_n$ on $B(p, 4r)$ satisfy $f_j(p) = 0, j=1,\cdots, n$ and
$\int_{B_r(p)} f_j\bar f_k = 0$ for $j\neq k$,
then there exists $1 \leq l\leq n$ so that
$$\frac{\fint_{B_{2r}(p)}|f_l|^2}{\fint_{B_r(p)}|f_l|^2}\geq
2^{10n\epsilon+2}.$$
\end{lemma}
\smallskip
\begin{proof} Now we prove Proposition \ref{STOS}.
We have to show that if $p\in W$ is smooth in complex analytic sense, then it is smooth in the metric sense.
Now  suppose that  there are holomorphic coordinates $(z_1,\cdots, z_n)$ on $B_{r_0}(p)$ with $z_j(p)=0, j=1,\cdots,n$ for some small constant $r_0$. Let us also fix a sequence of Gromov-Hausdorff approximation $B_i\to B_{r_0}$. By Lemma  \ref{Approximatinglemma}, on each $B_i$, we find $n$ holomorphic functions $z^i_j,j=1,\cdots,n$ converging uniformly to $z_j,j=1,\cdots,n$ (shrink the ball $B_{r_0}$ if necessary). Now we consider functions $v_i:=|dz^i\wedge d\bar z^i|_{g_i}:=|dz^i_1\cdots\wedge dz^i_n\wedge d\bar z^i_1\cdots\wedge d\bar z^i_n|_{g_i}$. Notice that 

$$v_i=\frac{dz^i\wedge d\bar z^i}{\Omega_{X,i}\wedge \bar\Omega_{X,i}}|\Omega_i\wedge \bar\Omega_i|_{g_i},$$
where $\Omega_{X,i}:=r_i^{-n}\Omega_X$ and $\Omega_X$ is a holomorphic volume form on $X$, whose existence  is guaranteed by our definition of Calabi-Yau variety in Definition \ref{a}. Here we assume $K_X$ is trivial and  the case $qK_X$ is trivial can be proved similarly.
By the Ricci flat condition, we may assume $|\Omega_i\wedge \bar\Omega_i|_{g_i}=1$. Therefore $v_i=|f_i|^2$, where $f_i=\frac{dz_1^i\cdots\wedge dz_n^i}{\Omega_{X,i}}$ is a holomorphic function on the regular part of $B_i$. By normality of $X$, $f_i$ can be extended to whole $B_i$. By the gradient estimate, $v_i=|f_i|^2$ is uniformly bounded, so $f_i$ converge to a holomorphic function $f$ on $B_{r_0}$. On the other hand, on the metric smooth part of $B_{r_0}$, we have smooth convergence of complex structure and metric, so $v_i$ converge to $|dz\wedge d\bar z|_{g_{C(Y)}}$ which is equal to $|f|^2$ on the metric regular part of $B_{r_0}$. Since $|dz\wedge d\bar z|_{g_{C(Y)}}\neq 0$ on the metric regular part and the metric singular part has Hausdorff codimension at least 4, $f$ as a holomorphic function, nonvanish on $B_{r_0}$. Otherwise, the zero of $f$ will be of complex codimension $1$.

 Now suppose $p$ is not a metric regular point, we will follow the argument of \cite[Proposition 4.1]{LiuSze1} to show that $f(p)=0$, which contradicts the nonvanishing of $f$ on $B_{r_0}$. By  rescaling and orthogonalization of $z_j,j=1,\cdots,n$ , we may assume 
 $$\int_{B_{r_0}}z_j\bar z_k=0, j\neq k,$$
 $$z_j(p)=0,\fint_{B_{r_0}}|z_j|^2=1, j=1,\cdots, n.$$  By rescaling the metric, we may assume without loss of generality that $r_0 = 2.$ By the gradient estimate $|dz\wedge d\bar z|_{g_{C(Y)}}\leq C$, where $C$ is the constant in Lemma \ref{gradholo}. Let $E:=\mathbb C\{z_1,\cdots, z_n\}$ be the complex vector space spanned by $z_i$ . On  $E$, there are two norms, given by
$L^2$ integration over $B_1(p)$ and $B_2(p).$ By a simultaneous diagonalization with
respect to the above two norms, we assume that $z_j, j=1,\cdots, n$ are also  orthogonal on $B_1(p).$
Fix the constant $\epsilon$ in Lemma \ref{criteria}. 

We claim that there is a constant $r_1$ such that if $r\leq r_1$, then for any non zero holomorphic function $g$ on $B_2(p)$ satisfying $g(p)=0$, we have
\begin{equation}\label{A}
\frac{\fint_{B_{2r}(p)}|g|^2 }{\fint_{B_{r}(p)}|g|^2 }
\geq 2^{2-\epsilon} 
\end{equation}
If the claim fails, then there are a sequence of positive number $r_i\to 0$ and non zero holomorphic function $g_i$ vanishing at $p$ such that 
$$\frac{\fint_{B_{2r_i}(p)}|g_i|^2 }{\fint_{B_{r_i}(p)}|g_i|^2 }
< 2^{2-\epsilon}. $$  By  rescaling, we may assume that  $\fint_{B_{r_i}(p)}|g_i|^2=1$ (Rescaling does not affect the above inequality). Then the rescaled metric space $(B_{2r_i}(p),\frac{d}{r_i})$ will converge to the radius $2$ ball of a metric cone $(V,o)$  and the corresponding sequence of functions $g_i|_{B_{2r_i}(p)}$ will converge to a non constant sublinear growth harmonic function on the metric cone $V$, which is a contradiction. Without lose of generality, we may assume $r_1>1$.
Now we apply Lemma \ref{criteria} to $z_j,j=1,\cdots, n$, then there exists $1\leq l\leq n$ such that
\begin{equation}\label{B}
\frac{\fint_{B_{2}(p)}|z_l|^2 }{\fint_{B_{1}(p)}|z_l|^2 }
\geq 2^{10n\epsilon+2} 
\end{equation}

Define $w_j=z_j2^{-\epsilon}$ for $j\neq l$ and $w_l=z_l2^{(n-1)\epsilon}$. So
$$dw_1\wedge\cdots\wedge dw_n=dz_1\wedge\cdots\wedge dz_n.$$
Moreover, from (\ref{A}), (\ref{B}) we have
$$\fint_{B_{1}(p)}|w_j|^2 < 2^{-\epsilon}, $$
for all $j$. By the gradient estimate we obtain that on $B_{\frac{1}{2}}(p)$
$$|dw_1 \wedge \cdots \wedge dw_n|=
|dz_1 \wedge \cdots \wedge dz_n|\leq
 C2^{-0.5n\epsilon}\leq C2^{-\epsilon}.$$ Notice that by (\ref{A}), we have  $\fint_{B_{2^{-i}}(p)}|z_j|^2\leq 2^{(\epsilon-2)i}$  for any positive integer $i$,  so by iteration, we obtain that for all $0 < r < 1$,
$|dz_1\wedge\cdots\wedge dz_n| \leq 2Cr^\epsilon$ on $B_r(p)$. Recall that $|f|^2=|dz_1\wedge\cdots\wedge dz_n|$ is nonvanishing. Contradiction. Proposition \ref{STOS} is proved.
\end{proof}
\begin{lemma}\label{klt} Any tangent cone
$C(Y)$ is homeomorphic to a normal affine variety. Moreover, $C(Y)$ has klt singularity. 
\end{lemma}
\begin{proof}
The normality is proved in Lemma \ref{normal} and affine property is a consequence of the cone structure (cf. \cite[Theorem 1.1]{LiSze2}). By general theory (cf. \cite[Definition 5.1 and 5.3]{EGZ}) and Proposition \ref{STOS} (The metric smooth part coincides with the complex analytic smooth part) , to show that $C(Y)$ has klt singularity, it suffices to show that for any point $p$, there is an open neighbourhood $U$ of $p$ and also an integrable nowhere zero holomorphic volume form over the complex analytic regular part $U_{reg}$ of $U$. 
Indeed, we will show that there is holomorphic volume form $\Theta$ on $C(Y)$ satisfying $\omega^n_{C(Y)}=\Theta\wedge\bar\Theta$, where $\omega^n_{C(Y)}=\ddbar r^2$ and $r$ is the distance function to the vertex of tangent cone. Recall that by the Ricci flat equation, on $X$, we have $\omega^n_\infty=\Omega_X\wedge\bar\Omega_X$ by choosing suitable rescaling of $\Omega_X$, where $\Omega_X$ is a holomorphic volume form on $X$. 
We can also rescale the equation as $r_i^{-2n}\omega^n_\infty= r_i^{-n}\Omega_X\wedge r_i^{-n}\bar\Omega_X$. Notice that on the smooth part of $C(Y)$,  the complex structure and also the metric $r_i^{-2}g_\infty$ converge smoothly, when $r_i\to 0$. Therefore, the covariant constant holomorphic volume form $r_i^{-n}\Omega_X$ on the metric space $(X,r_i^{-2}g_\infty)$ will converge smoothly to a holomorphic volume $\Theta$ defined on the smooth part of $C(Y)$. The  Ricci flat equation $\omega^n_{C(Y)}=\Theta\wedge\bar\Theta$ is preserved. The lemma is proved.
\end{proof}
\bigskip

\subsection{Tangent cone is affine variety and Rigidity of holomorphic spectrum of tangent cone} 

In this subsection, we discuss the holomorphic spectrum of tangent cones.

It follows from Remark \ref{T} that any tangent cone is indeed isomorphic to one of the element in $\mathcal C_p$. By Liu-Szekelyhidi \cite[Proposition 4.4]{LiSze2},  the Reeb vector field $\nu:=Jr\partial r$ generates a parameter group of isometry $\sigma_t$ acting on $C(Y)$. Taking  closure of the $R$ action $\sigma_t$ in the isometry group of $C(Y)$ preserving $o$, there is a torus $\mathbb T$ action on the tangent cone $(C(Y),o)$, which also acts holomorphically on the regular  part of the tangent cone.  For any $d>0$, let $\mathcal H_d$ be the vector space
of polynomial growth holomorphic functions $f$ on $C(Y)$ of degree at most $d$, i.e.,
satisfying $|f(x)|< Cd(o,x)$ for a constant $C>0$.  $\mathcal H_d$ is finite dimensional (cf. \cite{LiSze2}) and it decomposes into weight spaces under the
$\mathbb T$-action i.e. $f\in\mathcal H_d$ can be decomposed into a sum of eigenfunctions $f=\sum f_{\alpha_1}+\cdots+f_{\alpha_n}$ for $\alpha_k\in Lie(\mathbb T)^*$, where $e^{it}.f=e^{i\langle \alpha, t \rangle} f$ for $t\in Lie(\mathbb T)$. In particular,  $\nu.f_{\alpha}=\langle\xi,\alpha\rangle f_\alpha$, where $\xi\in Lie(\mathbb T)$ is induced by the Reeb vector field $\nu$.

Now take  a $\mathbb{T}$-invariant neighborhood of $o$ and denote it by $\Omega$. 
We can define degree of $f\in\mathcal O(\Omega)$ as $d(f):=min\{\langle\alpha,\xi\rangle|\alpha\in Lie(\mathbb T)^*, f_\alpha\neq 0\}$. There is also another way to characterize the degree, based on a three circle type lemma (cf. \cite{DS2}) as follows:
 \begin{equation}\label{degree2}
d(f)=\lim_{r\rightarrow0} \frac{\log \sup_{B_r(p)}|f(x)|}{(\log r)}.
\end{equation}


  Given  a tangent cone $C(Y)\in \mathcal C_p$, denote $R(C(Y))$ to be the ring of polynomial growth holomorphic function on $R(C(Y))$. Now we define the so called holomorphic spectrum $\mathcal S=\{d|d=d(f)\,\,\, \textnormal{for some} f\in\mathcal R(C(Y))\}$ and $\mu_d=\dim\{f\in R(C(Y))|d(f)=d\}$ for $d\in\mathcal S$.
A important fact is that the holomorphic spectrum of tangent cones is rigid. 

\begin{theorem}\cite[Theorem 3.3]{DS2}\label{RigidS} 
The holomorphic spectrum $\mathcal S:=\mathcal S(C(Y))$ and the Hilbert function of $C(Y)$ are independent of the  tangent cones in $\mathcal C_p$. 
 \end{theorem}
\begin{proof}
We only outline the proof. First of all, by \cite[Proposition 2.21]{DS2}, the holomorphic spectrum $\mathcal S$ for any tangent cone is contained in the set of algebraic numbers based on the result of \cite{GMSY}. Fix a tangent cone  $C(Y)\in\mathcal C_p$, by using Lemma \ref{Approximatinglemma} and Lemma \ref{L2Converge}, it holds that for any positive number $D\notin\mathcal S(C(Y ))$, there is a
neighborhood $U$ of $C(Y)$, which also depends on $D$, so that for all tangent cone $C(Y')$ in $U$, we have $\dim E_{D}(C(Y))=\dim E_{D}(C(Y'))$, where $E_{D}(C(Y)):=\bigoplus_{0<d< D}R_d(C(Y))$. Using the fact that $\mathcal C_p$ is compact and connected, there is a dense subset $\mathcal I$ of $\mathbb R$ such that, for any $D\in\mathcal I$, $N_D:=\dim E_D(C(Y))$ is invariant for all tangent cones. Then for each fixed $D\in\mathcal I$, we can define map $\Phi:\mathcal C_p\to R^{N_D}$ by sending a tangent cone to the degrees of holomorphic functions contained in $E_D$. $\Phi$ is continuous and its image is contained in a discrete set. So the image must be a point, which finish the proof.
\end{proof}
\bigskip
\subsection{Vanishing order of holomorphic functions by using the spectrum of tangent cone} In this subsection, we will define the degree of a holomorphic function defined on a germ of $X_\infty$ near $p$, hence giving a filtration of the local ring at $p$.

We remark that up to here, we have established all the necessary analytic estimates as in \cite{DS2} and the rest of the proof will be the same as \cite{DS2}. We now sketch the proof for reader's convenience. We will focus on how the uniqueness of tangent cone can be proved and skip the proof of several versions of three circle type lemmas.

We first fix some notations.
In this subsection, denote the pointed space  $(X_\infty, g_i=\lambda^{-2i}g_\infty,p_i=p)$ by $X_i$ for simplicity and denote  the unit ball centered at $p_i$ in $X_i$  by $B_i$. Let $\Lambda_i:B_i\to B_{i-1}$ be the inclusion map. Now Given a function $f$ defined on $B_{i-1}$, we also denote by $\Lambda_i. f$ the induced function on $B_i$ by restricting $f$ to $B_i$. 
If a sequence of $X_i$ (with different metric structure) converges by sequence to a tangent cone $C(Y)$, then $\Lambda_i$ converges to the dilation map $\Lambda: C(Y)\to C(Y)$ correspondingly. 


 

Let $\O_p$ the local ring of holomorphic functions defined in a neighborhood of $p\in X_\infty$.
Let us recall some notations from \cite{DS2}. Define  $||f||_{i}:=\int_{B_i}|f|^2$ and define $||[f]||_i$ to be the normalization of $f_i$ such that $||[f]_i||_i=1$. By the convergence of $B_i$ to $B$, holomorphic function $f\in\mathcal O_p$ will satisfy the following three circle type estimates, which holds on the tangent cone.


\begin{lemma}\label{3Circle}
For any given $\bar d\notin \mathcal S$, we can find $\mathbb N_0$ depending on $\bar d$ such that for all $i>j\geq \mathbb N_0$ and any non-zero holomorphic function $f$ defined on $B_j$,  if $\|f\|_{j+1}\geq \lambda^{\bar d} \|f\|_{j}$, then 
$\|f\|_{i+1}>\lambda^{\bar d} \|f\|_{i}$. 
\end{lemma}
\begin{proof} The proof is identical to the proof of \cite{DS2}. As remarked a couple of times before, the only difference is that we are using a sequence of singular K\"ahler-Einstein ball $B_i$ to approximate the unit ball $B$ on the tangent cone. On the other hand, by Lemma \ref{Approximatinglemma} and Lemma \ref{L2Converge}, things still work. We sketch the proof here.
 Suppose the conclusion is not true, then we can find a subsequence $\{\alpha\}\subset \{i\}$, and non-zero  holomorphic functions $f_\alpha$ defined on $B_{\alpha}$
with 
\begin{align*}
\|f_\alpha\|_{\alpha+1}\geq& \lambda^{\bar d} \|f_\alpha\|_{\alpha},\\
 \|f_\alpha\|_{\alpha+1}\geq& \lambda^{-\bar d}\|f_\alpha\|_{\alpha+2}.
\end{align*}
 By taking a subsequence we may assume $B_{\alpha}$ converges to a unit ball $B$ in some tangent cone. We may assume $||f_\alpha||_{\alpha+1}=1$. Then $||f_\alpha||_{\alpha}\leq \lambda^{-\bar d}$. Using the gradient estimate Lemma \ref{gradholo}, we may assume $f_{\alpha}$ converges  uniformly to a  function $F$ on $B$ by taking a further subsequence (Since $f$  is always well defined on a ball strictly larger than the unit ball $B_\alpha$). Hence we have
 \begin{align*}
 \|F\|_{L^2(B)}\leq \lambda^{-\bar d}, \|\Lambda. F\|_{L^2(B)}=1, \|\Lambda^2. F\|_{L^2(B)}=\lim_{\alpha\rightarrow\infty}\|f_{\alpha}\|_{\alpha+2}\geq \lambda^{\bar d}. \end{align*}
  Then by Lemma \cite[Lemma 4.4]{DS2},  $F$ must be  a homogeneous holomorphic function with degree $\bar d$ on $B$. This contradicts $\bar d\notin\mathcal S$.  

\end{proof}

The above three circle type estimate  implies the following consequence, which says that there is well defined notion of degree of $f\in\mathcal O_p$ (cf. \cite[Proposition 3.7]{DS2}).

Given $f\in\mathcal O_p$ and $f\neq 0$, then
the limit 
\begin{equation}\label{fdegree}d(f):=\lim_{i\rightarrow\infty}\frac{\log( \|f\|_{i+1}/\|f\|_{i})}{ (\log\lambda)}
\end{equation}
exists. Moreover $d(f)\in \mathcal S$ and $[f]_i$ converges strongly  by sequence to a non-zero  holomorphic functions, with degree $d(f)$, on the tangent cones. 

For $D\notin\mathcal S$ large, an orthonormal basis of homogeneous elements in $$E_{D}(C(Y)):=\bigoplus_{0<d< D}R_d(C(Y))$$ with give an equivalent embedding $\Phi$ of the tangent cone to some affine variety in $\mathbb C^N$. 
In equation (\ref{fdegree}), we see that there is a well-defined notion of degree of function $f\in\mathcal O_p$. We shall use homogeneous elements of $\mathcal O_p$ to construct a sequence of holomorphic map $F_\alpha$ from $B_\alpha$ to $\mathbb C^N$ converging to $\Phi$ when $B_\alpha$ converge to some  tangent cone $C(Y)$. ($\alpha$ is an increasing subsequence of $\mathbb Z$.) This motivates the following definition of adapted sequence in \cite{DS2}.
\begin{definition}\label{adaptedS}\textbf{(Adapted Sequence)} Let $P$ be a  finite dimensional subspace of $\mathcal O_p$. Let $n=\dim P$. 
An adapted sequence of bases consists of a basis $\{G_i^1, \cdots, G_i^n\}$ of $P$ for all large $i$,  with the following holds
\begin{enumerate}
\item For all $1\leq a\leq n$, $\|G_i^a\|_i=1$ and for $a \neq b$, $\lim_{i\rightarrow \infty} \int_{B_i} G_i^{a}\overline{G_i^{b}}=0$,
\item  $\Lambda_i. G^a_{i-1}=\mu_{ia} G^a_i+p_{i}^a$ for some $\mu_{ia}\in \mathbb C$, and $p_i^a$  $\in\mathbb C\langle G_i^1, \cdots, G_i^{a-1}\rangle$, with $\|p_i^a\|_i\rightarrow 0$,
\item  $\mu_{ia}\rightarrow \lambda^{d_a}$ for some $\lambda^{d_a}\in\mathcal S$.  
\end{enumerate}

\end{definition}

Suppose $P$ is a subspace of $\mathcal O_p$ with an adapted sequence of bases and  $f\in\mathcal O_p$ but not in $ P$. We want to construct an adapted sequence on the space spanned by $P$ and $f$. Let $\Pi_j(f)$ denotes the $L^2$ orthogonal projection of $f|_{B_j}$ to the orthogonal
complement of $P|_{B_j}$. Now we define   modulo $P$ degree as:
\begin{align*}
d_P(f):=\lim_{j\rightarrow\infty} \frac{\log(\|\Pi_{j+1} f\|_{j+1}/\|\Pi_i f\|_{j})}{(\log\lambda)}.\end{align*} This is similar to the definition (\ref{fdegree}). Moreover, we have the following lemma for modulo $P$ degree.
\begin{lemma}\cite[Proposition 3.12]{DS2} \label{Fadapted} Suppose $P$ has an adapted sequence of basis and $f\in\mathcal O_p$ is a nonzero function not included in $P$, then
 $d_P(f)\in \mathcal S$ and $P\oplus \mathbb C\langle f\rangle$ also has an adapted sequence of bases. 
    \end{lemma}

\subsection{Local tangent cones}
In this subsection, we study the local tangent cone and show that it degenerates to tangent cone in the Hilbert Scheme. We write the spectrum $\mathcal S$ of tangent cone (independent of the choice tangent cone) in the increasing order as $0=d_0<d_1<\cdots$ and also let $\mu_k=dim R_{d_k}(C(Y))$.


  For each $d\in \mathcal S$, we let $I_k$ be the subspace of $\mathcal O_p$ consisting of functions with $d(f)\geq d_k$.  Then there is a natural filtration

\begin{equation}\label{filtration}
\O_p=I_{0}\supset I_{1}\supset I_{2}\supset\cdots.
\end{equation}
By equation (\ref{degree2}), it is easy to verify that $I_jI_k\subset I_l$ if $d_l\leq d_j+d_k$.  

Hence this defines a  multiplicative filtration of ideals on $\O_p$ and this filtration induces the following graded ring:

$$R_p=\bigoplus_{k\geq 0} I_{k}/I_{{k+1}}.$$
Using Lemma \ref{Fadapted}, we can find an adapted sequence of bases of $I_k/I_{k+1}$ for all $k>0$. The following proposition is crucial.
\begin{proposition} \label{AdaptedDecom}
For all $k\geq 0$, we can find a decomposition $I_{k}=I_{{k+1}}\oplus J_k$, such that $\dim J_k=\mu_{k}$, and $J_k$ admits an adapted sequence of bases with $d(J_k)=\{d_k\}$.
\end{proposition}
\begin{proof} The proof is essentially due to \cite{DS2}. Due to its importance, we sketch the proof. We prove it by induction on $k$. When $k=0$, $J_0$ is the space of constant functions. 
Now we assume the conclusion holds for all $j\leq k-1$.  Define a set $$\mathbb S:=\{J\subset I_k|J\cap I_{k+1}=\emptyset\textnormal \,\,\,{and} \,\,\,J\,\,\, \textnormal{admit an adapted sequence of bases}\}.$$   Obviously, $\dim J\leq \mu_k$ for all $J\in \mathbb S$, otherwise the adapted sequence will converge by sequence to orthogonal holomorphic functions with degree $d_k$ on some tangent cone $C(Y)$ and this will contradict that $u_k=dim R_{d_k}$. Now let $J_k$ be a maximal element in $\mathbb S$.

\textbf{Claim:} $\dim J_k=\mu_k$. If not, fix a subsequence $\{\alpha\}\subset\{i\}$ and the correspongding tangent cone $C(Y)$, then the adapted sequence of basis of $J_k$ will converge to  a  orthonormal set  of homogeneous holomorphic functions $G^1, \cdots, G^p$ on $B \subset C(Y)$ of degree $d_k$, with $p<\mu_k$. So we may fine a non zero function $f$ in $R_{d_k}(C(Y))$ orthogonal to the complex vector space $\mathbb C\{ G^1, \cdots, G^p\}$.
Then by Lemma \ref{Approximatinglemma}, for $\alpha$ large, there is a sequence of holomorphic functions $f_\alpha$ defined on $B_\alpha$ converging uniformly to $f$ as $\alpha\rightarrow\infty$. 
Define $P:=\bigoplus_{i\leq k} J_i$.
By  \cite[Proposition 3.7, 3.11]{DS2}, we have for $\alpha$ sufficiently large,  $d(f_\alpha)\leq d_k$ and $d_P(f_\alpha)\leq d_k$.   On the other hand, by Lemma \ref{Fadapted}, for $\alpha_0$ large, if we let $g:=f_{\alpha_0}$, we obtain an adapted sequence of bases on $\hat P=P\bigoplus \C\langle g\rangle$ with $d(\hat P)=d(P)\cup \{d_P(g)\}$.  
Therefore $d_{P}(g)$ must be $d_k$. Otherwise,   an adapted  sequence of bases on $P\bigoplus\mathbb\langle g\rangle $ will converge by sequence to an orthogonal set of holomorphic functions with degree no bigger than $d_k$ on some tangent cone $C(Y)$. But this contradicts  the induction hypothesis that $\dim J_j=\mu_j$ for all $j\leq k-1$.  The same argument with the fact that $f$ is orthogonal to $\mathbb C\langle G^1,\cdots,G^p\rangle$, we will also have $d(g)=d_k$. Then $J_k\cup\{g\}$ is an  element in set $\mathbb S$ larger than $J_k$. Contradiction.

With the claim in hand,  we will finish the proof by showing that $I_{{k}}=I_{k+1}\bigoplus J_k$. Given any $f\in I_k$ and $f\notin J_k$, then by Lemma  \ref{Fadapted} we obtain a sequence of adapted bases on $J_k \bigoplus\mathbb C\langle f\rangle$. Using the fact that $\dim J_k=\mu_k$, we know $d_{J_k}(f)>d_k$. This implies $f\in I_{k+1}\bigoplus J_k$, and hence completes the proof.
\end{proof}
Suppose for some subsequence $\{\alpha\}\subset \{i\}$, $(X_\infty,\lambda^{-2\alpha}g_\infty)\to C(Y)$ in the Gromov-Hausdorff sense.
By previous discussion, we know $R(C(Y))$ is generated by $E_{D}(C(Y))$ for some large $D$.  An orthonormal basis of $E_D(C(Y))$ defines an equivariant embedding  $\Phi: C(Y)\rightarrow\mathbb C^N$.  Let $G_\xi$ be the group of linear transformations of $\mathbb C^N$ that commute with the $\mathbb T$ action, and let $K_\xi=G_\xi\cap U(N)$. 

Let  $P=\bigoplus_{0<k<D}J_k$.  By Proposition \ref{AdaptedDecom}, for $i$ large, we have an adapted sequence of bases of $P$ on $B_\alpha$, which in turn defines a sequence of holomorphic maps $F_\alpha: B_\alpha\rightarrow \mathbb C^N$.  converging uniformly to $\Phi$ when restricted to a ball (up to the $K_\xi$ action). By Lemma \ref{Fgeneric}, we may assume $F_\alpha$ is generically one-to-one for all $\alpha$ large. Using this property of $F_\alpha$, we shall prove the graded ring $R_p$ is finitely generated.

As in \cite{DS2}, we define two spaces: 
\begin{align*}S_k:&=\{ \textnormal{homogeneous polynomials on} \,\,\mathbb C^N \textnormal{with weighted degree}\,\, d_k\},  \\ V_k:&=\{ \textnormal{ kernel of the restriction map}\,\, S_k\rightarrow R_{d_k}(C(Y))\}.  \end{align*}Choose a splitting $S_k=V_k\bigoplus Q_k$, then we can identify $Q_k$ with $R_{d_k}(C(Y))$.   If we denote  the subspace of $\O_p$, consisting of those functions which are  pulled back from functions in $Q_k$ through map  $F_\alpha$, by $A_{k,\alpha}$.  By the generic one to one property of $F_{\alpha}$ , for $\alpha$ sufficiently large, we have if $f\in Q_k$
is non zero, then $F_\alpha^*f$ is non zero.  So  $A_{k,\alpha}$ is a subset of $I_k\setminus I_{k+1}$. Hence $dim A_{k,\alpha}=\mu_k$ and $I_k=I_{k+1}\bigoplus A_{k,\alpha}$. 
\begin{lemma}\cite[Lemma 3.15]{DS2}
$R_p:=\bigoplus_{k\geq 0}I_k/I_{k+1}$ is finitely generated by $\bigoplus_{k<D}I_k/I_{k+1}$.
\end{lemma}
 Since the graded ring $R_p$ is finitely generated, we my define the local tangent cone $W:=Spec (R_p)$ as an affine variety.

\subsection{An algebraic geometric argument to show the uniqueness.}

By proposition \ref{AdaptedDecom}, $R_p$ has the same grading as ring $R(C(Y))$, hence  $W$ admits an  action of $\mathbb T$.  Crucially, the chosen adapted sequence of bases of $P$ over $B_i$ will define a sequence of  embeddings of $W$ into $\mathbb C^N$, and we call the image $W_i$. As remarked in \cite{DS2},  there is a multi-graded
Hilbert scheme \textbf{Hilb}, which is a projective scheme parameterizing polarized
affine schemes in $\mathbb C^N$ invariant under the $\mathbb T$ action and with
fixed Hilbert function determined by ${\mu_k}$. Therefore, for i larege, $W_i$ and $C(Y)$ define points $[W_i]$ and $[C(Y)]$ in \textbf{Hilb} and all  $[W_i]$ are in the same $G_\xi$ orbit. Recall that if for some subsequence $B_\alpha$ converging to $B$,  the holomorphic map $F_\alpha$ induced by the adapted sequence of bases converge to $\Phi$, which is an embedding of the tangent cone to $\mathbb C^N$. We now show that $[W_\alpha]$ convergence to $C(Y)$ in the Hilbert Scheme.

\begin{lemma}\cite[Proposition 3.16]{DS2}\label{DEG} For $\alpha$ large, $W_\alpha$ is normal and
$[W_\alpha]$  converges to $[C(Y)]$ in \textbf{Hilb}, up to $K_\xi$ action.
 \end{lemma} 
 
\begin{proof}
We sketch the proof by  arguing as  \cite{DS2} for reader's convenience. Given an element $f\in V_k$, then for $\alpha$ large, we have  $F^*_\alpha f=g_\alpha+F^*_\alpha h_\alpha$ with $g_\alpha\in I_{k+1}$ and $h_\alpha \in  Q_k$. It is shown in \cite{DS2} that when $\alpha\to \infty$, $|h_\alpha|\to 0$ (here $|\cdot|$ can be an arbitrary norm on $S_k$). This is essentially because $f$ vanishes on the tangent cone. Now let $f_\alpha:=f-h_\alpha\in S_k$, then $F_\alpha^*f_\alpha=g_\alpha\in I_{k+1}$. So $f_\alpha$ vanishes on $W_\alpha$ (note that $F^*_\alpha f_\alpha$ lies in  $I_{k+1}$ while $f_\alpha$ in $S_k$). Now fix a constant $k_0$  such
that any ideal of $\mathbb C[x_1, \cdots, x_N]$ defining an element in \textbf{Hilb} is generated
by the homogeneous pieces of degree at most $k_0$. Then we for each function $g$ in  a basis of $V_k$ for all $k\leq k_0$, we may find a function $g_\alpha$ vanishing on $W_\alpha$ for $\alpha$ large. Notice that $g_\alpha\to g$, then it follows that $[W_\alpha]$  converges to $[C(Y)]$ in
  \textbf{Hilb}.
\end{proof}

 Recall that for all $C(Y')\in \mathcal C_p$,  a choice of orthonormal basis of $E_D(C(Y'))$ determines a holomorphic map $\Phi': C(Y')\rightarrow\mathbb C^N$ (not necessary an embedding). If $\Phi$ is an embedding for $C(Y)$, then for tangent cone $C(Y')$ close to $C(Y)$, then $\Phi'$ is indeed also an embedding by the following lemma.

\begin{lemma}\cite[Lemma 3.17]{DS2}\label{DEG2}
 There is a neighborhood $\mathcal U$ of $C(Y)$ in $\mathcal C_p$ such that for all $C(Y')\in \U$, $\Phi'(C(Y'))$ is normal.
 \end{lemma}

At last we prove our main results. In Lemma \ref{DEG},  we have constructed the degeneration of $[W_\alpha]$ to $[C(Y)]$ in the Hilbert Scheme and also in Lemma \ref{DEG2}, the embedding of nearby tangent cone. We can follow the algebraic geometric argument of 
\cite[cf. Page 353]{DS2} to complete the proof. For simplicity, we argue as  Li-Wang-Xu \cite{LWX}.
\smallskip

\textbf{Proof  of Theorem \ref{main1}:} We have proved that the tangent cone $C(Y)$ is Ricci flat K\"ahler cone with log terminal singularity, therefore it is K-polystable (cf. \cite[Corollary A.4]{LWX}). 
 By the equivalent degeneration of $[W_\alpha]$ to $[C(Y)]$ in \textbf{Hilb},  the proof of \cite[Theorem 1.4, page 61]{LX} shows that  that $W$ is indeed a K-semistable log Fano cone (We refer the reader to \cite{LX} for the precise definition of K-semistable log Fano cone).  Then by \cite[Theorem 1.2]{LWX}, the K-semistable log Fano cone $W$ can only degenerate to a unique K-polystable log Fano cone in the algebraic geometry sense. Therefore the metric tangent cone must be unique as a polarized affine variety. Then by the  uniqueness of Ricci flat K\"ahler cone metric (cf. \cite[Proposition 4.8]{DS2}), the tangent cone is also unique up to isometry.

\section{uniqueness of tangent cone: canonical polarized case}
In this section we prove our second main theorem.
\begin{theorem}[=Theorem 1.2]\label{main2}
 Let $X$ be a canonical polarized  variety with crepant singularity and let $\omega_{KE}\in K_X$ be the negative K\"ahler-Einstein metric constructed in \cite{EGZ} with locally bounded potential, then for any point $p\in X$, the tangent cone of metric $\omega_{KE}$ at $p$ is unique.
\end{theorem}

To prove this we, we need a lemma for the regularity of the tangent cone, which is similar to Lemma \ref{codim4}, in the context of canonical polarized variety with crepant singularity. We should remark that it is much harder than the Calabi-Yau case, since the K\"aher-Einstein metric space can only be approximated by smooth K\"ahler manifold with Ricci bounded below. So Anderson and Cheeger-Naber type argument do not work here. Instead, we use the deep structure result of Tian-Wang for the so called almost K\"ahler-Einstein manifold.

\begin{definition} \label{AKE} A sequence of pointed closed K\"ahler manifold  $(X_j,g_j, x_j), j\in \mathbb Z^+ $ is called    an almost K\"ahler-Einstein sequence if it satisfies the following: 

\begin{enumerate}

\item $Ric(g_j) \geq -\lambda_jg_j$, with $0<\lambda_j\leq 1$. 
\smallskip

\item There exists $r_0, \kappa >0$ such that for all $j=1, 2, ...$, 
%
$$Vol(B_{g_j}(x_j, r_0) ) \geq \kappa,$$ 

\item  Let $g_j(t)$ be the solution of the normalized K\"ahler-Ricci flow 
$$\ddt{g_j(t)} = -Ric(g_j(t)) - \lambda_jg_j(t), ~ g_j(0) = g_j, $$ then 
$$ \lim_{j\rightarrow \infty} \int_0^1 \int_{X_j} |R(g_j(t)) + n\lambda_j| dV_{g_j(t)} dt =0. $$

\end{enumerate}

\end{definition}

We remark that in the introduction of Tian-Wang \cite{TW}, in the definition of almost Einstein sequence, all $\lambda_j$ are equal to $1$. Here we allow $\lambda_j$ to be a sequence of  bounded constants. 
The following deep regularity result for Gromov-Hausdorff limit of almost K\"ahler-Einstein manifolds is proved in \cite{TW}. (We notice that although it is stated under the assumption that $\lambda_j=1$, its proof indeed allows $\lambda_j$ to be a bounded sequence.).
\begin{proposition}\cite[Main Theorem 2]{TW} \label{TW1} Let $(X_j, g_j, x_j)$ be a sequence of almost K\"ahler-Einstein manifolds (cf.  Definition \ref{AKE}). Then $(X_j, g_j,x_j)$ converges by sequence to  a metric length space $(X_\infty, p_\infty, d_\infty)$ satisfying 

\begin{enumerate}

\item $\R$, the regular set of $X_\infty$, is a smooth open dense convex set in $X_\infty$, 

\smallskip

\item the limit metric $d_\infty$ induces a  smooth K\"ahler-Einstein metric $g_{KE}$ on $\R$ satisfying $Ric(g_{KE}  ) = -\lambda g_{KE} $, where $\lambda=\lim_{j\to\infty}\lambda_j,$  

\smallskip

\item the singular set $\mathcal{S}$ has Hausdorff dimension no greater than $2n-4$.

\end{enumerate}

\end{proposition}

Crucially, it is shown by Song \cite{S} that the K\"ahler-Einstein metric $\omega_{KE}$   on the singular canonical polarized variety $X$ in Theorem \ref{main2} and its rescaled space  can be realized as the Gromov-Hausdorff limit of a sequence of almost K\"aher-Einstein manifolds. To be compatible with the notation used in the proof of Theorem \ref{main1}, we denote the K\"ahler-Einstein metric space $(X,\omega_{KE})$ by $(X_\infty, g_\infty)$ and also denote the rescaled metric space by  $(X_\infty, r_i^{-2}g_\infty)$, where $r_i>0$ is a constant.
\begin{lemma}\label{AKEGOOD}\cite[Lemma 4.2]{S}
There exists a sequence of almost K\"ahler-Einstein manifolds $(X'_i, g'_i, p_i)$ converging to $(X_\infty,g_\infty,p)$ in the Gromov- Hausdorff sense. Moreover, any rescaled space  $(X_\infty, r_i^{-2}g_\infty, p)$ is also the Gromov-Hausdoff limit of a sequence of almost K\"ahler-Einstein manifolds.
 \end{lemma}

Now we prove two regularity results for any tangent cone $C(Y)$ of $(X_\infty,g_\infty)$ at a given point $p\in X_\infty$.
\begin{lemma}\label{codim4II}
 Recall that we define $(X_i,g_i,p_i):=(X_\infty, r^{-2}_ig_\infty, p)$. If  $$(X_i,g_i, p_i)\stackrel{GH}{\longrightarrow}(C(Y),g_{C(Y)},o)$$
 in the Gromov-Hausdorff sense, then $g_i$ converge smoothly to the metric smooth part of $C(Y)$ and the metric singular set of $C(Y)$ has Hausdorff codimension at least 4.
\end{lemma}
\begin{proof}By a diagonal argument, this is direct consequence of Lemma \ref{AKEGOOD} and Proposition \ref{TW1}.
\end{proof}
 
We also show that the metric singular point of $X_i$ could not converge to the metric smooth part of $C(Y)$, which is an analogue of Lemma \ref{snotr}. In fact, we show a little bit more.
\begin{lemma}\label{snotrII}
Suppose we are in the context of Lemma \ref{codim4II}. There exists a dimensional constant $\eta(n)>0$ such that if $q_i$ is a sequence of points in the metric singular part of $X_i$ converging to a point $q$ in $C(Y)$, then at $q$, we have
$$\lim_{r\to0}\frac{Vol(B_r(q))}{Vol(B_r^{\mathbb E}(0))}<1-\eta(n).$$Moreover, if $q'\in C(Y)$ is point with $\lim_{r\to0}\frac{Vol(B_r(q'))}{Vol(B_r^{\mathbb E}(0))}>1-\eta(n),$ then $q'$ is smooth point in the metric sense.
\end{lemma}
\begin{proof}
Suppose the conclusion fails. For the given $\eta$ (a dimensional constant to be determined later), there exists $r_0$ such that 
$$\frac{Vol(B_{r_0}(q))}{Vol(B^{\mathbb E}_{r_0}(0))}>1-\eta,$$ where $\eta$ is a dimensional constant to be determined. Then by Colding's volume continuity Theorem, we have 
$$\frac{Vol(B_{r_0}(q_i))}{Vol(B^{\mathbb E}_{r_0}(0))}>1-2\eta$$ for $i$ sufficiently large. Here $B_{r_0}(q_i)$ is the raidus $r_0$ ball under the metric $r_i^{-2}g_\infty$. Using Lemma \ref{AKEGOOD}, for each fixed $i$, we can  fix a sequence of almost K\"ahler-Einstein manifolds  $(X'_j, g'_j, p_j),j\in\mathbb Z^+$  converging to $(X_i,r_i^{-2}g_\infty, q_i)$ in the Gromov-Hausdorff sense.   Notice that $Ric(g'_j)>-r_ig_j$ for all $j$ and $r_i\to 0$, we may assume that $$\frac{Vol(B_{r_0}(q_i))}{Vol(B^{r_i}_{r_0}(0))}>1-3\eta$$
for $i$ sufficiently large, where $B^{r_i}_{r_0}(0)$ is the radius $r_0$ ball in the complete simply connected manifold with Ricci curvature $-r_i$. Then by monotonicity of volume ratio, we have $$\frac{Vol(B_{s}(q_i))}{Vol(B^{r_i}_{s}(0))}>1-4\eta$$
for all $s<r_0$ and $i$ sufficiently large. So we have 
\begin{equation}\label{gap}\frac{Vol(B_{s}(q_i))}{Vol(B^{\mathbb E}_{s}(0))}>1-5\eta
\end{equation}for all $s<r_0$ and $i$ sufficiently large.



Notice that $X_i$ is the Gromov-Hausdorff limit of a sequence of almost K\"ahler-Einstein manifolds.
Now we need a gap theorem.  Since $q_i$ is a metric singular point of $X_i$ by our assumption,  by applying \cite[Theorem 2, item (4)]{TW} to the space $(X_i,r_i^{-2}g_\infty, q_i)$, we conclude that if $C(Z)$ is a tangent cone of $X_i$ at point $q_i$,   then
$$d_{GH}((B_1(q_i), g_{C(Z)}),(B_1(0), g_{\mathbb E})) > \epsilon(n),$$
where $\epsilon(n)$ is a dimensional constant. Simple contradiction argument implies that, there is a constant $r_1$ depending on $X_i$ such that if $s<r_1$, then $$d_{GH}((B_1(q_i), g_{s^{-2}g_i}),(B_1(0), g_{\mathbb E})) > \epsilon(n),$$
where $g_i=r_i^{-2}g_\infty$. Using Cheeger-Colding's volume cone implies metric cone Theorem, this in turn implies that there exist constants $r_2$ depending on $X_i$,  such that if $s<r_2$, then 
\begin{equation}\label{gap1}\frac{Vol(B_{s}(q_i))}{Vol(B^{\mathbb E}_{s}(0))}<1-\eta(n),
\end{equation}
where $\eta(n)$ is a dimensional constant. Now if choose $\eta=\frac{\eta(n)}{5}$, then inequality (\ref{gap}) will contradict inequality (\ref{gap1}). The first part of the Lemma is proved.

It remains to prove that if $q'\in C(Y)$ with $\lim_{r\to0}\frac{Vol(B_r(q'))}{Vol(B_r^{\mathbb E}(0))}>1-\eta(n)$, then $q'$ is a metric smooth point. Indeed, by the first part of the proof, $q'$ must be a limit point of a sequence of smooth points $q'_i$ in $X_i$. Then by volume continuity, we may assume that all $q'_i$  have large volume density, then the proof follows from Anderson's gap Theorem for Einstein manifolds by  choosing another smaller constant $\eta(n)$ if necessary. 
\end{proof}
\textbf{Proof of Theorem \ref{main2}:}
With Lemma \ref{codim4II} and Lemma \ref{snotrII} in hand, we are able to achieve the $H$- condition on the tangent cone and also graft the data on tangent cone to the approximation sequence.  We still need three uniform $L^2$ type estimates. More precisely, we point out that the uniform $L^2$ estimates for sections of $mK_X, m\in\mathbb Z^+$, and the uniform  $L^2$ estimates for gradient of the sections of $mK_X$ are proved by in
\cite[Proposition 4.4]{S}. Also Hormander type $L^2$ estimate for solving $\bar\partial$ equation on the singular canonical polarized variety with crepant singularity is proved in \cite[Proposition 4.5]{S}. Now we have all the ingredients to derive a similar partial  $C^0$ type result as in Proposition \ref{C0}. The rest of the proof is completely the same as the proof of Theorem \ref{main1}.



\bigskip

\noindent {\textbf{Acknowledgements:} We would like to thank Professor Jian Song for suggesting this problem and many useful discussions. We thank Gabor Szekelyhidi for many useful communications and comments. At last, we also thank Professor Chi Li, Bin Guo, Song Sun  for answering questions.}

\end{document}